\documentclass[12pt,reqno]{amsart}

\usepackage{hyperref}
\usepackage{amsmath,amsfonts,amssymb}
\usepackage[latin1]{inputenc}
\usepackage{verbatim}

\newtheorem{theorem}{Theorem}[section]
\newtheorem{lemma}[theorem]{Lemma}
\newtheorem{prop}[theorem]{Proposition}
\newtheorem{cor}[theorem]{Corollary}
\newtheorem{rmk}[theorem]{Remark}
\newtheorem{defn}[theorem]{Definition}

\newcounter{defn}

\newcommand{\sect}{\vspace{3mm} \setcounter{equation}{0} \setcounter{defn}{0} \section}

\renewcommand{\Re}{{\rm Re }\;}
\renewcommand{\Im}{{\rm Im }\;}

\newcommand{\w}[1]{\langle {#1} \rangle}
\newcommand{\pf}{\noindent {\bf Proof. \hspace{2mm}}}
\newcommand{\ef}{ \hfill $ \Box $ \vskip 3mm}
\newcommand{\be}{\begin{equation}}
\newcommand{\ee}{\end{equation}}
\newcommand{\bea}{\begin{eqnarray}}
\newcommand{\eea}{\end{eqnarray}}
\newcommand{\beas}{\begin{eqnarray*}}
\newcommand{\eeas}{\end{eqnarray*}}
\newcommand{\ep}{{\epsilon}}
\newcommand{\f}{\frac}
\newcommand{\na}{\nabla}

\newcommand{\bC}{{\mathbb C}}
\newcommand{\bR}{{\mathbb R}}
\newcommand{\bN}{{\mathbb N}}

\newcommand{\vB}{{\mathcal B}}

\newcommand{\vF}{{\mathcal F}}

\newcommand{\vH}{{\mathcal H}}

\newcommand{\vL}{{\mathcal L}}

\newcommand{\vN}{{\mathcal N}}

\newcommand{\vS}{{\mathcal S}}

\newcommand{\gM}{{\mathfrak M}}

\newcommand{\gm}{{\mathfrak m}}

\newcommand{\sgn}{\mbox{\rm sgn }}

\def\p{\partial}

\def\f{\frac}
\def\Dl{\Delta}
\def\na{\nabla}
\def\la{\lambda}

\def\vp{\varphi}

\def\Dl{\Delta}

\def\i{\infty}
\def\s{\sqrt}

\newcommand{\beq}{\begin{equation}}
\newcommand{\eeq}{\end{equation}}
\newcommand{\beqs}{\begin{equation*}}
\newcommand{\eeqs}{\end{equation*}}
\newcommand{\bal}{\begin{aligned}}
\newcommand{\eal}{\end{aligned}}

\begin{document}

\title[The Kramers-Fokker-Planck equation]{On the Kramers-Fokker-Planck equation \\
 with decreasing potentials in dimension one}
\author{Radek Novak \and Xue Ping Wang}

\address{R. Novak: Laboratoire de Mathématiques Jean Leray\\
Université de Nantes \\
44322 Nantes Cedex 3  France 
and
Department of Theoretical Physics\\
 Nuclear Physics Institute ASCR \\ 
 25068 \v{R}e\v{z}, Czech Republic.
 E-mail: r.novak@ujf.cas.cz }

\address{X.P. Wang: Laboratoire de Mathématiques Jean Leray\\
Université de Nantes \\
44322 Nantes Cedex 3  France.
E-mail: xue-ping.wang@univ-nantes.fr}


\begin{abstract} For quickly decreasing potentials with one position variable, the first threshold zero is always a resonance of the Kramers-Fokker-Planck operator. In this article we study low-energy spectral properties of the operator and calculate  large time asymptotics of solutions in terms of the Maxwellian.
\end{abstract}

\maketitle

\sect{Introduction}

The Kramers equation is a special Fokker-Planck equation describing the Brownian motion in an external field. This equation was derived and used  by H. A. Kramers \cite{kr} to describe  kinetics of chemical reaction. Later on it turned out that it had more general applicability to different fields such as supersonic conductors, Josephson tunnelling  junction and relaxation of dipoles (\cite{risc}). 
Mathematical analysis of the Kramers-Fokker-Planck (KFP, in short)  equation is initially motivated by  trend to  equilibrium for confining potentials (\cite{dv,hln,v}). Spectral problems of the KFP operator reveal to be quite interesting, because  this operator is neither elliptic nor selfadjoint.  After appropriate normalisation of physical constants and a change of unknowns,  the KFP  equation can be written into the form  
\be\label{equation} \p_t u(t; x,v)+P u(t; x,v)=0,\ (x,v)\in\bR^n\times\bR^n,   t >0,
\ee
with initial data
\be\label{initial} u(0; x,v)=u_0(x,v). 
\ee
Here $x$ and $v$ represent respectively position and velocity of the particle, $P$ is the KFP  operator given by
\be\label{operator} 
P=-\Dl_v+\f{1}{4}|v|^2-\f{n}{2} + v\cdot\na_x-\na V(x)\cdot\na_v, 
\ee
where the potential $V(x)$ is supposed to be  a real-valued $C^1$ function verifying
\be \label{ass1}
|V(x)| + \w{x}|\nabla V(x)| \le C\w{x}^{-\rho},  \quad x\in \bR^n,
\ee
for some $\rho \in \bR$ and  $\w{x} = (1 + |x|^2)^{1/2}$. 
Let $\gm$ be the function defined by
\be
\gm(x,v)= \f{1}{(2\pi)^{\f n 4}}e^{-\f 1 2 (\f{v^2}{2} + V(x))}.
\ee
Then $\gM=\gm^2$ is the Maxwellian (\cite{risc}) and $\gm$ verifies the stationary KFP equation
\be
P\gm =0 \quad  \mbox{ on } \bR^{2n}_{x,v}.
\ee
From the point of view of spectral analysis, large time behavior of solutions of (\ref{equation}) is closely related to low-energy spectral properties of $P$. 
If $V(x) \ge C|x|$ for some constant  $C>0$ outside some compact, then $\gm \in L^2(\bR^{2n}_{x,v})$ and
zero is a discrete eigenvalue of $P$. This case has been studied by many authors. It is known that
\be \label{eq1.6}
u(t) -\w{u_0, \gm} \gm = O(e^{-\sigma t}), \quad t\to +\infty,
\ee
in $L^2(\bR^{2n})$, where $\sigma >0$ can be evaluated in terms of spectral gap between zero eigenvalue and the real part of the other eigenvalues of $P$ and $V(x)$ is  normalized by
\[
\int_{\bR^n} e^{-V(x)} dx =1.
\]
See \cite{dv,hln,hrn,v} and references quoted therein. If $V(x)$ increases slowly: $V(x) \sim c \w{x}^\beta$ for some constants $c>0$ and  $\beta \in ]0, 1[$, then zero is an eigenvalue  embedded in the essential spectrum of $P$ and it is known that (\ref{eq1.6}) still holds with the right-hand side replaced by $O(t^{-\infty})$ (\cite{pc,dfg})  or more precisely  by $O(e^{-a t^{\f{\beta}{2-\beta}}})$ for some $a>0$ (\cite{lz}).
For decreasing potentials ( $\rho >0$ in (\ref{ass1}) ), zero is no longer an eigenvalue of $P$. It is proved  in \cite{w3} that for $n=3$ and  $\rho>2$,  one has 
\be \label{eq1.6b}
u(t)= \f{1}{(4\pi t)^{\f 3 2}}\w{u_0,\gm} \gm + O(\f 1{t^{\f 3 2 + \ep}}), \quad t\to +\infty, \ep >0, 
\ee
in some weighted spaces. 
(\ref{eq1.6b}) shows that for rapidly decreasing potentials, space distribution of particles is still governed by the Maxwellian, but the density of distribution decreases in times in the same rate as for heat propagation. Time-decay estimates of local energies are also obtained in \cite{w3} for  short-range potentials ($ \rho>1$) and in \cite{lz} for long-range potentials ($ 0 < \rho \le 1$). See also \cite{aggmms,epr,leb,lwx,n2} for other related works.
\\

In this work we study  one dimensional KFP equation with quickly decreasing potentials. 
It is known that for Schr\"odinger operators, low-energy spectral analysis in one and two dimensional cases is more difficult than higher dimensions and needs specific methods (\cite{abgh,bo,jn}) because zero is already a threshold resonance of the Laplacian in dimension one and two.  For the KFP operator with decreasing potentials,  the notion of thresholds and  threshold resonances is discussed in \cite{w3}. Although $\gm$ always verifies  the stationary KFP equation $P\gm =0$, a basic fact is that 
$\w{x}^{-s}\gm \not \in L^{2}(\bR^{2n})$ if  $n \ge 3$ and $1<s < \f n 2$, while $\w{x}^{-s}\gm \in L^{2}(\bR^{2n})$ for any $s>1$ if $n=1,2$.
In language of threshold spectral analysis, this means that for $n \ge 3$, zero is not a resonance of $P$ while for $n=1,2$, zero is a resonance of $P$ with $\gm$ as a resonant state. This is the main difference between the present work and \cite{w3}.
\\

Set $P = P_0 + W $ where 
\be
 P_0 = v\cdot\na_x-\Dl_v+\f{1}{4}|v|^2-\f{n}{2} \mbox{ and }  W= -\na_x V(x)\cdot\na_v.
\ee
$P_0$ and $P$ are  regarded as operators  in $L^2(\bR^{2n})$ with the maximal domain.  They are then maximally accretive.  Denote $e^{-tP_0}$ and $e^{-t P}$, $t \ge 0$, the strongly continuous semigroups generated  by $-P_0$ and $-P$, respectively.
If $\rho>-1$, $W$ is a relatively compact perturbation of the free KFP operator $P_0$: $ W (P_0+1)^{-1}$ is a compact operator in $L^2(\bR^{2n})$.
One can prove that
\be
\sigma_{\rm ess}(P)= \sigma(P_0) = [0, +\infty[
\ee 
and that non-zero complex eigenvalues of $P$ have positive real parts and may accumulate only towards points in $[0, +\infty[$. It is unknown for decreasing potentials whether or not the complex eigenvalues  does accumulate towards some point  in $[0, +\infty[$.
\\

 The main result of this work is the following
 
 \begin{theorem} \label{th1.1} Let $n=1$ and $\rho >4$. Then for any $s>\f 5 2$, there exists some $\ep>0$ such that
\be \label{11}
e^{-tP} = \f{1}{(4 \pi t)^{\f 1 2} }\left(\w{\cdot, \gm} \gm   + O(t^{-\ep})\right), \quad t \to +\infty 
\ee
as operators from $\vL^{2,s}$ to $\vL^{2,-s}$, where
\[
\vL^{2, r} = L^2(\bR^2_{x,v}; \w{x}^{2r} dxdv), \quad r\in \bR.
\]
\end{theorem}

To prove (\ref{11}),  the main task is to show that  the resolvent $R(z) = (P-z)^{-1}$ has an asymptotics of the form
\be \label{12}
R(z) =\f{i}{2 \s z} \w{\cdot, \gm} \gm  + O( |z|^{-\f 1 2 +\ep})
\ee
 as operators from $\vL^{2,s}$ to $\vL^{2,-s}$, for $z$ near zero and $z\not\in \bR_+$. Although (\ref{12}) and the decay assumption on the potential look the same as the  resolvent asymptotics of one dimensional Schr\"odinger operators in the case where zero is a resonance but not an eigenvalue (\cite{abgh,bo,jn}), its proof is quite different from  the Schr\"odinger case.  In fact, the known methods for the Schr\"odinger operator can not be applied to the KFP operator, mainly because the perturbation $W$ is a first order differential operator. In this work we use the method of \cite{w3}  to calculate the low energy  asymptotic expansion  for the free resolvent $R_0(z) =(P_0-z)^{-1}$ of the form
\be
R_0(z) = \f{1}{\s z} G_{-1} + G_0 + \s z G_1 + \cdots 
\ee
in appropriate spaces, where $G_{-1}$ is an operator of rank one. By a careful analysis of the space $\vN$ of resonant states of $P$ defined by (\ref{N}), we prove that $1+ G_0W$ is invertible on $\vL^{2,-s}$, $s> \f 3 2$. (\ref{12}) is derived from the equation
\be
R(z) = D(z) (1 + M(z))^{-1}R_0(z)
\ee
for $z$ near zero and $z\not\in \bR_+$, where 
\[
D(z) = (1+ R_1(z)W)^{-1} \mbox{ with  } R_1(z) = R_0(z) - \f{1}{\s z} G_{-1}
\]
 and 
\[
M(z) =\f{1}{\s z} G_{-1} W  D(z).
\]
As in threshold spectral analysis for Schr\"odinger operators, a non-trivial problem here is to compute  the value of some spectral constants involving the resonant state of $P$. Indeed, in most part of this work  only the condition $\rho>2$ is needed. The stronger assumption $\rho>4$ is used to show that some number $m(z)$ is nonzero for $z$ near $0$ and $z\not\in \bR_+$ (see (\ref{mz})), which allows to prove the invertibility of  $1+ M(z)$ and to calculate its inverse.\\

The organisation of this article is as follows. In Section 2, we recall some known results needed in this work. 
The low-energy asymptotics of the free resolvent in dimension one is calculated in Section 3. The threshold spectral analysis of $P$ is carried out in Section 4. We prove in particular that zero resonance is simple and $1+G_0W$ is invertible. The low-energy asymptotics of the full resolvent (\ref{12}) is proved in Section 5,  which implies in particular that if $\rho>4$, zero is not an accumulation point of complex eigenvalues of $P$. Finally,  Theorem \ref{th1.1} is  deduced in Section 6  by using  a high-energy resolvent estimate of \cite{w3} valid in all dimensions.
\\

{\bf \noindent  Notation.} For $r \ge 0$ and $s \in \bR$, introduce the weighted Sobolev space
\[
\vH^{r,s} = \{u \in \vS'(\bR^{2n}); (1 - \Delta_v + |v|^2 + \w{D_x}^{\f 2 3})^{\f r 2 } \w{x}^s u \in L^2\}.
\]
For $r <0$, $\vH^{r,s} $ is defined as the dual space of $\vH^{-r,-s} $ with the dual product identified with the scalar product of $L^2$. The natural norm on $\vH^{r,s}$ is denoted by $\|\cdot\|_{r,s}$.
When no confusion is possible, we use $\|\cdot\|$ to denote the usual norm of $L^2(\bR^{2n})$ or that bounded operators on $L^2$. Set $\vH^r = \vH^{r,0}$ and $\vL^{2,s} =\vH^{0,s}$.
Denote $\vB(r,s; r',s')$ the space of continuous linear operators from $\vH^{r,s}$ to $\vH^{r',s'}$. 
The weighted Sobolev spaces $\vH^{r,s}$ are introduced in accordance with the sub-ellipticity of $P_0$:  although $P_0$ does not map $\vH^{1,s}$ to $\vH^{-1,s}$,
the sub-elliptic estimate of $P_0$ (Corollary \ref{cor2.4}) implies that $(P_0 +1)^{-1} \in \vB(-1,0; 1, 0)$ and a commutator argument shows that $(P_0 +1)^{-1} \in \vB(-1,s; 1, s)$ for any $s\in \bR$. \\

\sect{Preliminaries}

In this Section we fix notation and state some known results which will be used in this work. Denote by $P_0$ the free KFP operator (with $\nabla V=0$):
\be
P_0=v\cdot\na_x-\Dl_v+\f{1}{4}|v|^2-\f{n}{2}, \quad x, v \in \bR^n.
\ee
In terms of Fourier transform in $x$-variables, we have for $u \in D(P_0)$
\bea
 (P_0 u) (x, v)& = &\vF_{x\rightarrow\xi}^{-1}\hat{P}_0(\xi) \hat{u}(\xi, v), \quad \mbox{ where }\\
\hat{P}_0(\xi)& = &-\Dl_v+\f{v^2}{4} -\f{n}{2}+ i v\cdot\xi, \\
 \hat{u}(\xi, v) & = & (\vF_{x\rightarrow\xi}u)(\xi, v) \triangleq \int_{\bR^n} e^{-i x\cdot\xi}u(x, v) \; dx.
 \eea
Denote
\be
D(\hat{P}_0) =  \{f \in L^2(\bR^{2n}_{\xi, v}); \hat{P}_0(\xi)f   \in L^2(\bR^{2n}_{\xi, v})\}.
\ee
Then $\hat{P}_0  \triangleq\vF_{x\rightarrow\xi} P_0 \vF_{x\rightarrow\xi}^{-1}$ is a direct integral of the  family of complex harmonic operators $\{\hat{P}_0(\xi); \xi \in \bR^n \}$.
 \\

For fixed $\xi \in \bR^n$,  $\hat{P}_0(\xi)$ can be written as
\[
\hat{P}_0(\xi) =-\Dl_v+\f{1}{4}\sum^n_{j=1}(v_j+2i\xi_j)^2-\f{n}{2}+|\xi|^2.
\]
 $\{\hat{P}_0(\xi), \xi\in \bR^n\}$ is  a holomorphic family of type A  with constant domain
 $D= D(-\Dl_v+\f{v^2}{4})$ in $L^2(\bR^n_v)$.  Its spectrum and eigenfunctions can be explicitly calculated.
 Let $F_j(s)=(-1)^je^{\f{s^2}{2}}\f{d^j}{ds^j}e^{-\f{s^2}{2}}, j \in \bN,$ be the Hermite polynomials and
$$
\varphi_j(s)=(j!\sqrt{2\pi})^{-\f{1}{2}}e^{-\f{s^2}{4}}F_j(s)
$$
the normalized Hermite functions.  For $\xi \in \bR^n$ and $\alpha=(\alpha_1, \alpha_2, \cdots, \alpha_n) \in \bN^n$, define
\be
\psi_\alpha(v) = \prod_{j=1}^n\varphi_{\alpha_j}(v_j) \mbox{ and } \psi_\alpha^\xi(v) = \psi_\alpha(v + 2i \xi).
\ee
One can check (\cite{w3}) that the spectrum of $\hat{P}_0(\xi)$ is given by
\be \label{e2.7}
\sigma (\hat{P}_0(\xi)) =\{ l +\xi^2; l \in \bN\}.
\ee
Each eigenvalue $l +\xi^2$ is semi-simple ({i.e.}, its algebraic multiplicity and geometric multiplicity are equal) with multiplicity $m_l = \# \{\alpha\in \bN^n; |\alpha| =\alpha_1 + \alpha_2 + \cdots + \alpha_n = l\}$.
The Riesz projection associated with the eigenvalue $l + \xi^2$ is given by
 \be\label{projection}
 \Pi^\xi_l\phi=\sum_{\alpha, |\alpha|=l} \langle  \phi, \psi^{-\xi}_\alpha  \rangle\psi^\xi_\alpha, \quad \phi \in L^2.
  \ee

 The following result is useful to study the boundary values of the resolvent $R_0(z) = (P_0 -z)^{-1}$.
Let $\hat R_0(z) = ( \hat P_0-z)^{-1}$ and $\hat R_0(z, \xi) = ( \hat P_0(\xi)-z)^{-1}$ for $z \not \in \bR_+$. Then
$R_0(z) = \vF_{x\to\xi}^{-1} \hat R_0(z) \vF_{x\to\xi}$.  

\begin{prop}\label{prop2.1}  Let $l\in \bN$ and $ l < a <l +1$ be fixed. Take  $\chi \ge 0$ and $\chi \in C_0^\infty(\bR^n_\xi)$ with supp $\chi \subset \{\xi, |\xi|\le a + 4\}$, $\chi(\xi) =1$ when $|\xi| \le a +3$ and $ 0 \le \chi(\xi) \le 1$.  Then one has
\be \label{R0}
\hat R_0(z, \xi) = \sum_{k=0}^l \chi(\xi) \frac{\Pi_k^\xi}{\xi^2+k -z} + r_l(z,\xi),
\ee
for any $\xi \in \bR^n$ and $z \in \bC$ with $\Re z <a$ and $\Im z \neq 0$. Here $r_l(z, \xi)$ is holomorphic in $z$ with $\Re z < a $ verifying the estimate
\be \label{reste}
\sup_{\Re z <    a, \xi \in \bR^n} \|r_l(z, \xi)\|_{\vL(L^2(\bR^n_v))} <\i.
\ee
\end{prop}

See Proposition 2.7 of \cite{w3} for the proof.  As a consequence of Proposition \ref{prop2.1} and known results for the boundary values of the resolvent of $-\Delta_x$, we obtain the following

\begin{cor}\label{cor2.2} Let $n \ge 1$ and  $R_0(z) = (P_0 -z)^{-1}$, $z \not \in \bR_+$. \\

(a). With the notation of  Proposition \ref{prop2.1}, one has
\be \label{RR0}
R_0(z) = \sum_{k=0}^l b_k^w(v, D_x, D_v)(-\Delta_x+k -z)^{-1} + r_l(z)
\ee
where $r_l(z)$ is $\vB(L^2)$-valued holomorphic function for $\Re z <a$ and $ b_k^w(v, D_x, D_v)$ is the Weyl pseudo-differential operator with symbol $b_k(x, \xi, \eta)$ given by
\be \label{Bk}
b_k(v, \xi,\eta) = \int_{\bR^n} e^{-i v'\cdot\eta/2}\left(\sum_{|\alpha|=k} \chi(\xi) \psi_\alpha( v+ v' + 2i\xi)\psi_\alpha( v- v' + 2i\xi) \right) dv'.
\ee
 In particular,
\be \label{B0}
b_0(v, \xi,\eta) = 2^{\f n2}e^{-v^2 -\eta^2 + 2iv\cdot\xi + 2\xi^2} \chi(\xi).
\ee
\\

(b).   Let $I$ be a compact interval of $\bR$ which does not contain any non negative integer. Then for any $s>\f 1 2$, one has
\be \label{LAP1}
\sup_{\lambda \in I; \epsilon \in ]0, 1]}\|R_0(\lambda \pm i\epsilon)\|_{\vB(-1,s; 1,-s)} <\infty
\ee
The boundary values of the resolvent $R_0(\lambda \pm i0) = \lim_{\ep \to 0_+} R_0(\lambda \pm i\epsilon)$ exist
in $\vB(0,s; 0,-s)$ for $\lambda \in I$ and is H\"older-continuous in $\lambda$.
\end{cor}

Seeing (\ref{RR0}), it is natural to define $\bN$ as set of thresholds of the KFP operator $P$ (\cite{w3}). 
Note that an exponential upper-bound in $\la$ for  $R_0(\lambda \pm i \ep)$, $\ep>0$ fixed, is obtained
in \cite{lz}  by method of harmonic analysis in Besov spaces. \\

For high energy resolvent estimate, we need the following result proved in Appendix A.2 of \cite{n2}.
 \\

\begin{theorem}\label{th2.3}  There exists some constant $C>0$ such that
\be \label{3.2.10}
\| (1 -\Delta_v + v^2 + |\xi|^{\f 2 3} + |\la|^{\f 1 2}) (\hat{P}_0(\xi) + \f n 2 + 1 -i \la)^{-1}\| \le C
\ee
uniformly in $\xi \in \bR^n$ and $\la \in \bR$.
\end{theorem}

As consequence, we obtain a uniform sub-elliptic estimate for the free KFP operator.

\begin{cor} \label{cor2.4}  One has
\be \label{3.3.1}
|\la| \|u\|^2 + \|\Delta_v u \|^2 + \||v|^2u\|^2 + \| |D_x|^{\f 23 } u\|^2 \le C  \|(P_0 + \f{n+2} 2 - i\la) u \|^2 , \quad  
\ee
for $u\in \vS(\bR_{x,v}^{2n}) $ and $\la \in \bR$.  In addition, $P_0$ defined on $\vS(\bR^{2n}_{x,v})$ is essentially maximally accretive.
\end{cor}

Let us indicate that the essential maximal accretivity of $P_0$ is discussed in \cite{n}. Henceforth we still denote by $P_0$ its 
closure in $L^2$ with maximal domain $D(P_0)= \{u \in L^2(\bR_{x,v}^{2n}); P_0u \in L^2(\bR_{x,v}^{2n})\}$. 
To determine the spectrum of $P_0$ which is unitarily equivalent with a direct integral of $\hat{P}_0(\xi)$, $\xi\in \bR^n$, in addition to (\ref{e2.7}), one needs a resolvent estimate uniform with respect to $\xi\in \bR^n$ proved
in \cite{w3}: $\forall z\in \bC\setminus \bR_+$,
\be \label{e2.9}
\sup_{\xi\in \bR^n} \| ( \hat P_0(\xi)-z)^{-1}\| \le C_z.
\ee
See \cite{d4} for  the necessity of such uniform resolvent estimate in order to determine the spectrum of direct integral of a family of non-selfadjoint operators. (\ref{e2.7}) and (\ref{e2.9}) show that
\be \label{e2.10}
\sigma(P_0) = \cup_{\xi \in \bR^n} \sigma(\hat{P}_0(\xi)) =[0, +\infty[.
\ee
Under the condition (\ref{ass1}) on $V$ for some  $\rho>-1$, $|\nabla V(x)| \to 0$ as $|x|\to +\infty$.
By Corollary \ref{cor2.4}, $W = - \nabla V(x)\cdot \nabla_v$ is  relatively compact with respect to $P_0$. It follows that
\be
\sigma_{\rm ess} (P) =[0, +\infty[
\ee
and discrete spectrum of $P$ is at most countable with possible accumulation points included in $[0, +\infty[$.
\\

\sect{The free resolvent in dimension one}
  
We use (\ref{RR0}) with $l=0$ to calculate the asymptotics of $R_0(z)$  near the first threshold zero. 

\begin{prop}\label{prop3.4.1} Let $n =1$. One has the following low-energy resolvent asymptotics for $R_0(z)$:
for $s, s'> \f 1 2$,
there exists $\ep>0$ such that
\be \label{eq3.4.1}
R_0(z) = \f{1}{\s z}(G_{-1} + O(|z|^\epsilon)), \mbox{ as }  z\to 0, z\not\in\bR_+,
\ee
as operators in $\vB(-1,s; 1, -s')$.  More generally, for any integer $N \ge 0$ and $s >N +\f 3 2$,  there exists $\epsilon>0$
\be \label{eq3.4.2}
R_0(z) =\sum_{j= -1}^N z^{\f j 2}G_j + O(|z|^{\f N 2 + \epsilon}), \mbox{ as }  z\to 0, z\not\in\bR_+,
\ee
as operators in $\vB(-1,s; 1, -s)$. Here the branch of $z^{\f 1 2}$ is chosen such that its imaginary part is positive when $z\not \in \bR_+$ and $G_j \in \vB(-1,s; 1, -s)$ for $s>j+\frac 3 2$, $j \ge 0$.  In particular,
\bea \label{3.4.3}
G_{-1} &= &\f i 2  \w{\cdot, \gm_0}\gm_0 \\
G_0 &= & F_0 + F_1, \label{3.4.4}
\eea
where 
\be
\gm_0(x, v) = 1\otimes \psi_0(v)
\ee
with $\psi_0(v) = \f{1}{(2\pi)^{\f 1 4}} e^{-\f{v^2} 4}$ the first eigenfunction of harmonic oscillator, $F_0$ is the operator with integral kernel
\be \label{3.4.5}
F_0(x,v; x', v') = -\f 1 2 \psi_0(v)\psi_0(v') |x-x'|
\ee
 and $F_1\in \vB(-1,s; 1,-s')$ for any $s, s' > \f 1 2$.
\end{prop}
\pf 
For $z \not \in \bR_+$, (\ref{RR0}) with $l=0$  shows that
\be \label{eq3.4.31}
 R_0(z) =  b_0^w(v,D_x,D_v)(-\Delta_x  -z)^{-1} + r_0(z),
\ee
with $ r_0(z) \in \vB(-1,0; 1,0)$  holomorphic in $z$  when  $\Re z < a$ for some $a \in ]0, 1[$. Here the cut-off $\chi(\xi)$ is chosen such that  $\chi\in C_0^\infty(\bR^n)$ and $\chi(\xi) = 1$ in a neighbourhood of $\{ |\xi|^2 \le a\}$. Therefore $r_0(z)$ admits a convergent expansion in powers of $z$ for $z$ near $0$
\[
r_0(z)   = r_0(0) +  z r_0'(0) + \cdots + z^n \f{r_0^{(n)}(0)}{n!} + \cdots
\]
in $\vB(-1,0; 1,0)$. It is sufficient to study the lower-energy expansion of $ b_0^w(v,D_x,D_v)(-\Delta_x  -z)^{-1}$.
\\

Note that in one dimensional case, the integral kernel of the resolvent $(-\Delta_x -z)^{-1}$ is given by
\be
\f{i}{2\s z} e^{i\s z |x-y|}, z\not\in \bR_+, x,y \in \bR
\ee
where the branch of $ \s z$ is chosen such that its imaginary part is positive for $z\not\in \bR_+$. 
The integral kernel of $  b_0^w(v,D_x,D_v)(-\Delta_x  -z)^{-1}$, $z \not\in \bR_+$,  is given by
\be \label{3.4.6}
K(x,x';v,v';z) =\f{i}{2\s z} \int_{\bR} e^{i\sqrt{z}|y-(x-x')|} \Phi(v,v',y) \; dy
\ee
with
\bea
\Phi(v,v',y) &= &(2\pi)^{-\f 3 2}e^{-\f 1 4 (v^2 + v'^2)} \int_{\bR} e^{i(y-v-v')\cdot \xi + 2\xi^2}\chi(\xi) \; d\xi \nonumber\\
 & = & \psi_0(v)\psi_0(v') \Psi(y-v-v') \label{Psi}
\eea
where $\Psi$ is the inverse Fourier transform of $e^{2\xi^2}\chi(\xi)$.
 Since $\chi\in C_0^\infty$, one has the following asymptotic expansion for $K(x,x';v,v';z)$ : for  any $\ep \in [0,1]$ and $N\ge 0$
\be \label{eq3.4.7}
|K(x,x';v,v';z) - \sum_{j=-1}^N z^{\f j 2} K_j(x, x', v, v') | \le C_{N, \ep} |z|^{\f{N+\ep}{2}}|x-x'|^{N+1+\ep}e^{-\f 1 4 (v^2 + v'^2)}
\ee
where
\bea
\label{3.4.8}
K_j(x, x';v,v')& = &  \frac{i^{j+2}}{2 \; (j+1)!} \int_{\bR}|y-(x-x')|^{j+1} \Phi(v,v',y) dy.
\eea
Remark  that for $N\ge 0$, $s', s > N + \f 1 2 $ and $0 <\ep < \min\{s,s'\}- N- \f 1 2$ and $\ep \in ]0, \f 1 2]$
\[
\w{x}^{-s}\w{x'}^{-s'}  |x-x'|^{N+\ep}e^{-\f 1 4 (v^2 + v'^2)} \in L^2(\bR^{4}).
\]
 We obtain the asymptotic expansion for
$ b_0^w(v,D_x,D_v)(-\Delta_x  -z)^{-1}$ in powers of $z^{\f 1 2}$ for $z$ near $0$ and $z\not\in \bR_+$.
\be \label{eq3.4. 34}
b_0^w(v,D_x,D_v)(-\Delta_x  -z)^{-1} =\sum_{j=-1}^N z^{\f j 2}K_j + O(|z|^{\f N 2 + \epsilon}), \mbox{ as }
\ee
as operators in $\vB(0,s'; 0, -s)$, $s', s > N + \f 3 2$. By the sub-elliptic  estimate of $P_0$, this expansion still holds in  $\vB(-1,s'; 1, -s)$. This proves (\ref{eq3.4.2}) with $G_k$ given by
\be
 G_{2j} = K_{2j} + \f{ r_0^{(j)}(0)}{j!}, \quad  G_{2j-1} = K_{2j-1}, \quad j \ge 0.
\ee

To show (\ref{3.4.3}) and (\ref{3.4.4}), note that since $\chi(0)=1$, one has
\[
 \int_{\bR} \Phi(v,v',y) \; dy =   \psi_0(v)\psi_0(v').
\]
The  first two terms in the expansion of $K(x,x';v,v';z)$ can be simplified as
\bea
K_{-1}(x, x',v,v') & = & \frac{i}{2} \int_{\bR}\Phi(v,v', y)  dy  
 =   \frac{i}{2}  \psi_0(v)\psi_0(v')  \\
K_0(x, x',v,v')&= & - \frac{1}{2} \int_{\bR} \Phi(v, v',y) |y-(x-x')|dy  \\
        &= &  - \frac{1}{2}  \psi_0(v)\psi_0(v') |x-x'| -\f 1 2  \int (|y-(x-x')|- |x-x'|) \Phi(v, v',y) dy.
         \nonumber
\eea
Therefore (\ref{3.4.3}) is true and $G_0$ can be decomposed as:  $G_0 = F_0 + F_1$  with $F_0$ defined by (\ref{3.4.5}) and $F_1 = K_{0,1} + r_0(0)$, $K_{0,1}$ being the operator with the integral kernel
\[
K_{0,1}(x,x', v, v') = -\f 1 2  \int_\bR (|y-(x-x')|- |x-x'|)  \Phi(v, v',y)   dy,
\]
which is a smooth function and
\[
K_{0,1}(x,x', v, v') = O(\psi_0(v)\psi_0(v'))
\]
for $|x-x'|$ large. Therefore $K_{0,1}$ is bounded in $B(-1,s; 1,-s')$ for any $s, s' > \f 1 2$.
This shows that $F_1 = K_{0,1} + r_0(0)$ has the same continuity property, which proves (\ref{3.4.4}).
\ef

\begin{cor} \label{cor3.4.1} Let $n=1$ and $e^{-t P_0}$, $t\ge 0$,  be the strongly continuous semigroup generated by $-P_0$. Then for any integer $N\ge 0$ and $s> 2N + \f 1 2$, the following asymptotic expansion holds for some $\ep>0$
\be \label{eq3.4.11}
e^{-tP_0} = \sum_{k =0}^N t^{-\f{2k+ 1}{ 2}} \beta_k G_{2k-1} + O(t^{- \f{2N+1} 2-\ep}), \quad t\to +\infty,
\ee
in $\vB(0,s, 0,s)$. Here $\beta_k$ is some non zero constant. In particular, the leading term $\beta_0G_{-1}$ is a rank-one operator given by
\be
\beta_0G_{-1} = \f{1}{(4 \pi)^{\f 1 2} }\w{\cdot, \gm_0} \gm_0 : \vL^{2,s} \to \vL^{2, -s}
\ee
for any $s >\f 1 2$. 
\end{cor}

The proof of Corollary \ref{cor3.4.1} uses Proposition \ref{prop3.4.1} and a representation formula of the semigroup $e^{-t P_0}$ as contour integral of the resolvent $R_0(z)$ in the right half-plane. See the proof of Theorem \ref{th1.1} for more details.\\

\sect{Threshold spectral properties}

 Assume that $V\in C^1(\bR^n; \bR)$ and 
\be
|V(x)| + \w{x} |\nabla V(x)|\le C \w{x}^{-\rho}
\ee
for some $\rho >0$.
Consider the null space of $P$ defined by
\be \label{N}
\vN =\{u; u \in \vH^{1, -s}, \forall s > 1 \mbox{ and }  Pu =0\}. 
\ee
 Since zero is not an eigenvalue of $P$, $\vN$ is the spaces of resonant  states of $P$ associated with zero resonance. See \cite{w3} for the definitions in general case. Remark that for $n=1$,  one can equally take $s>\f 1 2$ in the above definition, instead of $s>1$.  But the condition $s>1$ is necessary to define appropriately resonant states for $n=2$.  Clearly,  $\gm \in \vN$. We want to  prove that in one dimensional case, one has: $\dim \vN = 1$. In order to calculate the leading term of the resolvent expansion at threshold zero, we need also to calculate solutions of some integral equation.
\\

\begin{lemma} \label{lem3.4.3} Let $\rho>0$ and $n = 1$. 
If $u  \in \vH^{1, -s}$  for some $ s <\rho + \f 1 2 $ and satisfies  the equation $  Pu =0$, then 
\be
\w{Wu, \gm_0} =0,
\ee
where
\[
\gm_0(x, v) = 1\otimes \psi_0(v).
\]
\end{lemma}
\pf Suppose for the moment $n \ge 1$. Since $u \in \vH^{1, -s}$, one has $Wu \in \vH^{0, \rho+1-s} \subset L^2$. Using the equation $Pu=0$ and the ellipticity of $P$ in velocity variables $v$, we deduce that $(-\Delta_v + v^2)u(x,\cdot) \in L^2(\bR^n_v)$ a.e. in $x \in \bR^n$.
Taking scalar product of $Pu$ with $\psi_0(v)$ in $v$-variables, one has 
\[
 \w{(Pu)(x, \cdot), \psi_0}_v =  0, \quad\mbox{  a. e. }\quad x\in \bR^n.
\]
Since $\psi_0$ is the first eigenfunction of the harmonic oscillator in $v$, one has also
\[
\w{Pu, \psi_0}_v =  \w{v\cdot \nabla_x u, \psi_0}_v    - \w{\nabla_x V(x) \cdot \nabla_v u, \psi_0}_v
\]
{ a. e.} in $x \in \bR^n$. These two relations imply that
\be \label{4.11}
2 \nabla_x \cdot \w{ \nabla_v u, \psi_0}_v +  \nabla_x V(x)  \cdot\w{\nabla_v u, \psi_0}_v =0.
\ee
The above equation holds for $n \ge 1$. In the case  $n=1$, $\w{ \nabla_v u, \psi_0}_v$ is a scalar function in $x$ and the differential equation (\ref{4.11})
determines $\w{ \nabla_v u, \psi_0}_v$ up to some constant:
\be
\w{ \nabla_v u, \psi_0}_v = Ce^{-\f{V(x)}2}, \mbox{  a. e. in } x\in \bR
\ee
for some constant $C$. It is now clear that in one dimensional case, one has
\be
\w{Wu,\gm_0} = -\int_{\bR} V'(x) \w{ \partial_v u, \psi_0}_v  dx = -C\int_{\bR} V'(x) e^{-\f{V(x)}2} dx =0,
\ee
because $V(x) \to 0$ as $|x|\to +\infty$. \ef

Lemma \ref{lem3.4.3} is important in  threshold spectral analysis of the KFP operator in dimension one.
We believe that this result still holds when $n\ge 2$, but the last argument above does not hold if $n \ge 2$. In fact when $n \ge 2$, (\ref{4.11})
only implies that the vector-valued function $\w{ \nabla_v u, \psi_0}_v$ is of the form
\be
\w{ \nabla_v u, \psi_0}_v = e^{-\f{V(x)} 2} \overrightarrow{F}(x)
\ee
where $ \overrightarrow{F} \in L^2(\bR^n; \w{x}^{-2s} dx)$ and $\nabla \cdot \overrightarrow{F} =0$ in sense of distributions, which are not sufficient to conclude that $\w{Wu, \gm_0} =0$. \\

From now on, assume that $\rho>2$ and $n=1$.  Then by the sub-elliptic estimate for $P_0$, $G_0W$ is a compact operator in $\vH^{1, -s}$ for $ \f 3 2 < s< \f{\rho +1} 2$. We want to study  solutions of the integral equation
\be
( 1+ G_0W)u =\beta \gm_0
\ee
 for  $u\in \vH^{1,-s}$ and $\beta \in \bC$.\\

\begin{lemma}\label{lem3.4.4} Let $\rho>2$ and $u \in \vH^{1,-s}$ for some $ \f 3 2 < s <\f{\rho +1} 2$ such that  $ (1 + G_0W)u=\beta \gm_0$ for some $\beta \in \bC$. Then $Pu=0$. In particular, one has: $\w{Wu,\gm_0}=0$.
\end{lemma}
\pf One has seen that
\[
R_0(z) = \f{G_{-1}}{\sqrt{z}} + G_0 + o(1)
\]
in $\vB(-1,r; 1,-r)$ for any $r >\f 3 2$.  Therefore,
\[
 G_0 Wu = \lim_{z\to 0, z\not \in \bR_+} (R_0(z) - \f{G_{-1}}{\sqrt{z}})Wu
\]
in $H^{1,-r}$. Since $P_0G_{-1}=0$ in $\vH^{-1,r}$, one has for $\lambda <0$
\[
P_0 (R_0(\lambda) - \f{G_{-1}}{\sqrt{\lambda}})Wu = Wu +\lambda R_0(\lambda)Wu. 
\]
The m-accretivity of  $P_0$  implies
\[
\|\lambda R_0(\lambda)W\| \le 1, \quad \la <0.
\]
It follows that
\[
\|\lambda R_0(\lambda)Wu\| \le \|Wu\| \le C \|u\|_{1,-s},  \quad \f 3 2 < s < \f{ \rho +1} 2,
\]
uniformly in $\la <0$.  In addition,  if $\f 1 2 < s' <  \f{ \rho +1} 2$, one has
\[
\|\lambda R_0(\lambda)Wu\|_{1,-s'} \le  \| \la R_0(\lambda)\|_{\vB(0, s'; 0,-s')} \|W u \|_{0, s'} \le C|\la|^{\f 1 2} \|u\|_{1,-s}
\]
for $\la <0$. These two bounds show that
\be
w-\lim_{\lambda \to 0_-} \lambda R_0(\lambda)Wu =0, \quad \mbox{in } L^2(\bR^2).
\ee
Since  $u= -G_0Wu +\beta \gm_0$ and $P_0 \gm_0=0$,  the following equalities hold:
\[
P_0 u   = - w-\lim_{\lambda \to 0_-}  P_0 (R_0(\lambda) - \f{G_{-1}}{\sqrt{\lambda}})Wu  = -Wu
\]
in sense of distributions. This proves that $P u =0$. In particular Lemma \ref{lem3.4.3} shows that $\w{Wu,\gm_0}=0$. 
\ef

\begin{prop} \label{prop3.4.5} Let $u \in \vH^{1,-s}$ for some $ \f 3 2 < s <\f{\rho +1} 2$  such that  $ (1 + G_0W)u=\beta \gm_0$ for some $\beta \in \bC$. Then one has
\be
u(x,v) =(\beta - C_1(x) - v C_2(x))\psi_0(v) + r(x,v)
\ee
where $C_j \in L^\infty$ and  $C_j' \in L^1$, $j=1,2$,  and $(1+v^2- \partial^2_v) r\in L^2(\bR^2_{x,v})$. 
In addition, 
\be
\lim_{x\to \pm \infty} C_1(x) = \pm  d_1, \quad \lim_{x\to \pm \infty} C_2(x) = 0
\ee
where
\be
d_1= -\f 1 4 \int \int_{\bR^2 } ( x +\f {v}{ 2 } )\psi_0(v) \nabla V(x)\nabla_v u(x,v) dx dv.
\ee
In particular, $u \in \vH^{1,-s}$ for  any $s > \f 1 2$.
\end{prop}
\pf
Recall that $G_0 = K_1 + r_0(0)$ where $r_0(0)$ is bounded from $\vH^{-1}$ to $\vH^1$ and $K_1$ is an operator of integral kernel
\be
K_1(x, x'; v, v') = -\f 1 2 \int_\bR |y-(x-x')| \Phi(v,v';y) dy
\ee
with
\[
\Phi(v,v',y) = \f 1 2 \psi_0(v)\psi_0(v') \Psi( y-v-v'),
\]
$\Psi$ being the inverse Fourier transform of $e^{2\xi^2} \chi(\xi)$. Let $u \in \vH^{1,-s}$, $ \f 3 2  < s <\f{\rho +1} 2$,  such that  $ (1 + G_0W)u=\beta \gm_0$.
By Lemma \ref{lem3.4.4},
\be
\w{Wu, \gm_0} = 0.
\ee
Set $w= K_1Wu$. Then $u +w -\beta \gm_0 = - r_0(0)Wu$ belongs to $ L^2$. Let us study the asymptotic behavior of $w$ as $|x| \to \infty$.
Put
\[
F(x', y , v, v') = \psi_0(v)\psi_0(v') \Psi( y-v-v')\nabla V(x')\nabla_v u(x',v').
\]
Making use of the asymptotic expansion
\[
 |y-(x-x')| = |x-x'| -\f {y(x-x')}{  |x-x'|} + O( \f{y^2}{|x-x'|})
\]
for $|x-x'|$ large,  one obtains  that
\bea
w(x, v) &= & \f 1 4 \int\int_{\bR^3} |y-(x-x')| F(x', y , v, v') dy dx'dv' \nonumber \\
&\simeq & \f 1 4 \int \int_{\bR^3 } ( |x-x'| -\f {y(x-x')}{  |x-x'|} )F(x', y , v, v') dy dx'dv' \label{S1} \\
& = &  \f 1 4 \int \int_{\bR^2 } ( |x-x'| -\f {(v+v')(x-x')}{  |x-x'|} ) \psi_0(v)\psi_0(v') \nabla V(x')\nabla_v u(x',v') dx'dv'.
\nonumber
\eea
Here and in the following, ``$ \simeq$'' means the equality modulo some term in $L^2(\bR^2)$.
\\

Recall that since $\Psi$ is the inverse Fourier transform of $e^{2\xi^2} \chi(\xi)$, one has
\[
\int_\bR \Psi(y) dy =1, \quad \int_\bR y \Psi( y) dy = 0 
\]
and that according to  Lemma \ref{lem3.4.3}
\be
\int_{\bR^2} \psi_0(v') \nabla V(x')\nabla_v u(x',v') dx'dv' = -\w{Wu, \gm_0} = 0. \label{S2}
\ee
  The term related to $|x-x'|$ on the right-hand side of (\ref{S1}) is equal to
\bea
\lefteqn{\f 1 4 \int \int_{\bR^2 } ( |x-x'| )\psi_0(v)\psi_0(v') \nabla V(x')\nabla_v u(x',v') dx'dv'} \nonumber \\
&= & \f 1 4  \left(\int_{-\infty}^x - \int_{x}^{+\infty}\right) ( x-x' )\psi_0(v)\nabla V(x')\w{\nabla_{v'} u(x',\cdot), \psi_0}_{v'} dx'  \label{S3} 
\eea
Applying (\ref{S2}), one has for $x \le 0$
\bea
\lefteqn{\left|x\left(\int_{-\infty}^x - \int_{x}^{+\infty}\right ) \psi_0(v)\nabla V(x')\w{\nabla_{v'} u(x',\cdot), \psi_0}_{v'} dx'\right|}\nonumber  \\
& =& \left|2x\int_{-\infty}^x \psi_0(v)\nabla V(x')\w{\nabla_{v'} u(x',\cdot), \psi_0}_{v'} dx'\right | \nonumber  \\
& \le & C|x| \left\{\int_{-\infty}^x \w{x'}^{-2(\rho+1 -s)} dx' \right \}^{\f 1 2} \psi_0(v)\|u\|_{\vH^{1,-s}}\nonumber  \\
& \le & C' \w{x}^{-\rho +s + \f 1 2} \psi_0(v)\|u\|_{\vH^{1,-s}}
\eea
Since $\rho>2$ and $s< \f{\rho + 1} 2$ , this proves that the term
\[
x\left(\int_{-\infty}^x - \int_{x}^{+\infty}\right) \nabla V(x')\w{\nabla_{v'} u(x',\cdot), \psi_0}_{v'} dx' 
\]
is bounded for $x\le 0$ and tends to $0$ as $x \to -\infty$. The same conclusion also holds as $x \to +\infty$, using once more (\ref{S2}).
In the same way one can check that
\[
\left(\int_{-\infty}^x - \int_{x}^{+\infty}\right) x' \nabla V(x')\w{\nabla_{v'} u(x',\cdot), \psi_0}_{v'} dx'  
\]
is bounded for $x\in\bR$.
The other terms in (\ref{S1}) can be studied in a similar way. Finally we obtain that
\bea \label{S4}
w(x,v) &\simeq &   (C_1(x) + vC_2(x))  \psi_0(v) \mbox{ where } \\
C_1(x) &= & \f 1 4 \int \int_{\bR^2 } (x- x' -\f{v'} 2) \sgn(x-x')\psi_0(v') \nabla V(x')\nabla_v u(x',v') dx'dv' \label{S5}\\
C_2(x) &= & - \f 1 8 \int \int_{\bR^2 }\sgn (x-x') \psi_0(v') \nabla V(x')\nabla_v u(x',v') dx'dv'. \label{S6}
\eea
It follows from Dominated Convergence Theorem that the limits  
\be
\lim_{x\to \pm \infty} C_j(x) = \pm d_j
\ee
exist, where
\bea
d_1 & = & -\f 1 4 \int \int_{\bR^2 } ( x' +\f {v'}{ 2 } )\psi_0(v') \nabla V(x')\nabla_v u(x',v') dx'dv' \\
d_2  &= & -\f 1 8 \int \int_{\bR^2 } \psi_0(v') \nabla V(x')\nabla_v u(x',v') dx'dv' =0.
\eea
This proves that 
\[
 u \simeq  \beta \gm_0 -w \simeq (\beta - C_1(x) - vC_2(x)) \psi_0(v)
 \]
modulo some terms in $L^2(\bR^2)$.  In particular, $u \in \vH^{1,-s}$ for any $s>\f 1 2$.
Since $\rho>2$, one can also check that $C_j'(x)$ belongs to $L^1(\bR)$, $j=1,2$. \ef

\begin{theorem} \label{th3.4.6}  Assume $\rho >2$. If $u \in \vH^{1, -s}$, $\f 3 2 < s  < \f{\rho+1} 2$, satisfies the equation $ (1+ G_0 W)u =0$, then $u=0$.
\end{theorem}
\pf 
Let $\chi_1 \in C_0^\infty(\bR)$ be a cut-off with $\chi_1(\tau) =1$ for $ |\tau| \le 1$ and $\chi_1(\tau) =0$ for $|\tau| \ge 2$ and $0 \le \chi_1(\tau) \le 1$. Set $\chi_R(x) = \chi_1(\f x R)$ for $ R\ge 1$  and $u_R(x,v) = \chi_R(x) u(x, v)$.
Then one has
\[
P u_R = \frac{v}{R} \chi'(\frac{x}{R}) u.
\]
Taking the real part of the equality $\w{ P u_R, u_R}= \w{ \frac{v}{R} \chi'(\frac{x}{R}) u,  u_R}$, one obtains
\be  \label{eq4.6}
\int\int_{\bR^2} |(\partial_v + \frac v 2)u(x,v)|^2 \chi_R(x)^2 \; dx dv =  \w{ \frac{v}{R} \chi'(\frac{x}{R}) u, u_R}.
\ee
According to Proposition \ref{prop3.4.5}, $u$ can be decomposed as
\be
u(x,v) = z(x,v) + r(x,v)
\ee
where $z(x, v) = -(C_1(x) + vC_2(x)) \psi_0(v)$ and $C_1, C_2$ and $r$ are  given in Proposition \ref{prop3.4.5}.
Since $\psi_0(v)$  is even in $v$, the term $\w{ \frac{v}{R} \chi'(\frac{x}{R}) z, \chi_R z}$ is reduced to
\bea
\lefteqn{2 \Re \w{ \frac{v^2}{R} \chi'(\frac{x}{R})  C_1\psi_0,\chi_R C_2 \psi_0} }\\
& = &-\Re \int\int_{\bR^2} v^2 \psi_0(v)^2 \chi_R(x)^2 \f{d}{dx} ( C_1(x) \overline C_2(x)) dx dv  \\
&\to & -\Re \int\int_{\bR^2} v^2 \psi_0(v)^2  \f{d}{dx} ( C_1(x) \overline C_2(x)) dx dv =0 
\eea
as $R\to +\infty$, because  $\f{d}{dx} ( C_1(x) \overline C_2(x))$ belongs to $L^1$ and $C_1(x) \overline C_2(x)\to 0$ as $|x|\to +\infty$.
The term $|\w{ \frac{v}{R} \chi'(\frac{x}{R}) r, u_R}|$ can be estimated by
\[
|\w{ \frac{v}{R} \chi'(\frac{x}{R}) r,  u_R}| \le C R^{-(1-s)}\|u\|_{L^{2,-s}} \|\w{v}r\|_{L^2} 
\]
for $\f 1 2 < s <1$. Similar estimate also holds for $|\w{ \frac{v}{R} \chi'(\frac{x}{R}) z, \chi_R r}|$.
Summing up, we proved  that
\be
\lim_{R\to + \infty}   \w{ \frac{v}{R} \chi'(\frac{x}{R}) u, u_R} =0
\ee
which implies that  $(\partial_v + \frac v 2)u(x,v) =0$ a.e. in $x$ and $v$. Since $u \in \vH^{1, -s}$ for any $s>\f 1 2$ and $Pu=0$, 
it follows that $u$ is of the form $u(x, v) = D(x) e^{-\frac{v^2}{4}}$ for some $D \in L^{2,-s}(\bR)$ verifying  the equation
\be \label{eq4.71}
 D'(x) + \f 1 2   V'(x) D(x)=0
\ee
 in sense of distributions on $ \bR $.  It follows that $D(x) = \alpha  e^{-\f{V(x)}{2}}$ a.e. for some constant $\alpha$. Hence
\[
u(x, v) = \alpha e^{-\f{v^2}{4} -\f{V(x)}{2}}.
\]
In particular, one has
\bea
\int_0^R \int_{\bR_v} u(x,v) dv dx &= &\sqrt{\pi} \alpha R + O(1) \\
\int_{-R}^0 \int_{\bR_v} u(x,v) dv dx &= &\sqrt{\pi} \alpha R + O(1) 
\eea
as $R\to + \infty$. But according to  Proposition \ref{prop3.4.5}, one has for some constant $d_1$
\bea
\int_0^R \int_{\bR_v} u(x,v) dv dx &=& -\f{d_1}{\sqrt{2}} R + o(R)\\
\int_{-R}^0 \int_{\bR_v} u(x,v) dv dx &=& \f{d_1}{\sqrt{2}} R + o(R).
\eea
as $R\to +  \infty$. One concludes that $\alpha = d_1= 0$. Therefore $u=0$. 
\ef

Since $G_0W$ is a compact operator on $\vH^{1,-s}$, $\f 3 2 <s < \f{\rho+1} 2$,  it follows from Theorem \ref{th3.4.6} that
$1+G_0W$ is invertible and 
\be
(1+ G_0W)^{-1} \in \vB(1,-s; 1, -s).
\ee

\begin{theorem} \label{th3.4.7} Let $\rho>2$. One has:
\be \label{3.4.31}
\vN =\{ u \in \vH^{1,-s}; (1+G_0W) u =\beta \gm_0 \mbox{ for some } \beta \in \bC,  \f 3 2 < s< \f{\rho+1} 2 \}.
\ee
In particular, $\vN$ is of dimension one  and 
\be \label{3.4.32}
(1+G_0W) \gm = \gm_0
\ee
\end{theorem}
\pf  To prove (\ref{3.4.31}), it remains to  prove the inclusion 
\be \label{3.4.33}
\vN \subset \left \{ u \in \vH^{1,-s}; (1+G_0W) u =\beta \gm_0 \mbox{ for some } \beta \in \bC,  \f 3 2 < s< \f{\rho+1} 2 \right\}.
\ee
The inclusion in the opposite sense is a consequence of  Lemma \ref{lem3.4.4} and Proposition \ref{prop3.4.5}.\\

Let $u\in \vN$ and $\la <0$.  Then $u \in \vH^{1, -r}$ for $r>1$ and $r$ close to $1$ and $P_0 u = -Wu \in \vL^{2,\rho+1 -r}$. 
 By Corollary \ref{cor2.2}, the resolvent $R_0(\lambda)$ can be decomposed as
\be
R_0(\lambda) = b_0^w(v,D_x,D_v) (-\Delta_x- \lambda)^{-1} + r_0(\lambda)
\ee
where
\[
b_0(v, \xi, \eta) = 2^{\f 3 2}  e^{-v^2-\eta^2 + 2i v\cdot\xi + 2\xi^2} \chi(\xi)
\]
with $\chi$  a smooth cut-off around $0$ with compact support, and $r_0(\lambda)$ is uniformly bounded  as operators in $L^2$ for $\lambda < a$ for some $a\in ]0, 1[$.
 One has
\be \label{3.4.23}
u + R_0(\la)Wu = -\la R_0(\la)u = -\la \left( b_0^w(v,D_x,D_v) (-\Delta_x- \lambda)^{-1} + r_0(\lambda)\right )u
\ee
for $\lambda <0$.  Recall  the following estimate for $r_0(\lambda)$ (see (2.85) in \cite{w3} ):
\be \label{eq3.4.51}
\|\w{x}^{-s} r_0(\lambda) \w{x}^{s} f\| \le C (\|f\| + \|H_0 f\|)
\ee
for $f \in D(H_0)$, $\la <a$ and $s\in [0,2]$,  where $H_0 = -\Delta_v + v^2 -\Delta_x$.
It follows from \eqref{eq3.4.51} that
\be
\la r_0(\la)u = O(|\la|), \quad \la <0, 
\ee
in $\vH^{1,-r}$. \\

 Let $\phi \in \vS(\bR)$ such that $\int_\bR\phi(x) dx =1$.  Then 
\[
\Pi =\w{\cdot, \phi \otimes \psi_0}\gm_0
\]
is a projection  on $\vH^{1,-s}$ for any $s > \f 1 2$ onto the linear span of $\gm_0$. Set  $ \Pi' =1-\Pi$.
The term  $ \Pi'\lambda  b_0^w(v,D_x,D_v) (-\Delta_x- \lambda)^{-1} u$  can be evaluated as follows.
Making use of the inequality 
\[
|e^{-a} - e^{-b}| \le | a-b| (e^{-a} + e^{-b}), \quad a, b \ge 0,
\]
the quantity
\beas
\lefteqn{|\lambda  \Pi' b_0^w(v,D_x,D_v) (-\Delta_x- \lambda)^{-1} u (x,v)|}  \\
& =& 
\f{\sqrt{|\la|}} 2  \left|\int_{\bR^4}( e^{-\sqrt{|\lambda|}|y-(x-x')|} -  e^{-\sqrt{|\lambda|}|y-(y'-x')|}) \phi(y')\Phi(v,v',y)u(x',v')\; dy dy'dx' dv'  \right |
\eeas
is bounded by
\[
  |\la| \int_{\bR^4} |x-y'| (e^{-\sqrt{|\lambda|}|y-(x-x')|} +  e^{-\sqrt{|\lambda|}|y-(y'-x')|} )|\phi(y')\Phi(v,v',y)u(x',v')|\; dy dy'dx' dv'.
  \]
The integral involving the term $e^{-\sqrt{|\lambda|}|y-(x-x')|}$ can be evaluated as follows:
\beas 
\lefteqn{ |\la| \int_{\bR^4} |x-y'| e^{-\sqrt{|\lambda|}|y-(x-x')|} |\phi(y')\Phi(v,v',y)u(x',v')|\; dy dy'dx' dv'}
 \\ 
 & \le & C_1 (1+ |x|)|\la| \int_{\bR^3} e^{-\sqrt{|\lambda|}|y-(x-x')|} \big|\Phi(v,v',y)u(x',v')\big|\; dy dx' dv' \\
 & =  &  C_2 (1+ |x|)|\la| \int_{\bR^3} e^{-\sqrt{|\lambda|}|y-(x-x')|} \big|\psi_0(v)\psi_0(v') \Psi(y-v-v')u(x', v')\big|\; dy dx' dv'
 \\ 
 & \le & C_3 (1+ |x|)|\la|  \|u\|_{\vL^{2,-r}} 
 \\
  & & \times \left\{\int_{\bR^3} \big|\w{x'}^r e^{-\sqrt{|\lambda|}|y-(x-x')|} \psi_0(v)\psi_0(v') \Psi(y-v-v')\big|^2\; dy dx' dv'  \right\}^{\f 1 2} 
\\
& \le &  C_4 (1+ |x|)^{1+r}|\la|  \|u\|_{\vL^{2,-r}}  \left\{\int_{\bR} \big|\w{x'}^r e^{-\sqrt{|\lambda|}|x'|} \big|^2\;  dx' \right\}^{\f 1 2}\psi_0(v)
\\
& \le &  C_5 (1+ |x|)^{1+r}|\la|^{\f 3 4 - \f r 2}  \|u\|_{\vL^{2,-r}} \psi_0(v)
\eeas
for some constants $C_j$. A similar upper-bound also holds for the integral involving the term $ e^{-\sqrt{|\lambda|}|y-(y'-x')|}$.
Putting them together, we obtain a point-wise upper-bound
\be
 |\lambda  \left(\Pi' b_0^w(v,D_x,D_v) (-\Delta_x- \lambda)^{-1} u\right)(x,v)| \le C (1+ |x|)^{1+r}|\la|^{\f 3 4 - \f r 2} \psi_0(v) \|u\|_{\vL^{2,-r}} 
\ee
 This proves that for $ 1 < r < \f 3 2$, 
\be
 \lambda  \Pi' b_0^w(v,D_x,D_v) (-\Delta_x- \lambda)^{-1} u \to 0, \mbox{ as } \la \to 0_-
\ee
in $\vL^{2,-(\f 3 2 +r + \ep)}$, $\ep>0$.
Applying $\Pi'$ to \eqref{3.4.23} and taking the limit $\la \to 0_-$, we get
\be
\Pi'(1+G_0W)u =0.
\ee 
This means that there exists some constant $\beta\in \bC$ such that $(1+G_0W)u =\beta \gm_0$. The proof of  \eqref{3.4.31} is complete.
\\

Since $1+G_0W$ is injective, one deduces from \eqref{3.4.31} that $\vN$ is of dimension one. It is clear that $\gm \in \vN$ and \eqref{3.4.31}  implies that
\be
(1+G_0W)\gm =\beta \gm_0
\ee
for some $\beta\in \bC$. Proposition \ref{prop3.4.5} applied to $\gm$ shows that $\gm$ has asymptotic behavior
\[
\gm(x,v) = (\beta \mp d_1 + o(1))\psi_0(v), \quad x \to \pm \infty
\]
with $d_1\in \bC$ given in  Proposition \ref{prop3.4.5}. Comparing these relations with  the trivial expansion of  $\gm(x, v)$:
\[
\gm(x, v) = (1+ O(\w{x}^{-\rho}))\psi_0(v)
\]
 for $x \to \pm \infty$, one concludes that $\beta =1$ and $d_1=0$, which prove (\ref{3.4.32}).
\ef

\sect{Low-energy expansion of the resolvent}

Let $U_\delta =\{z; |z|<\delta, z\not\in \bR_+\}$, $\delta >0$,  and $\f 3 2 < s <\f {\rho+1} 2$. 
Recall that $(1+G_0W)^{-1}$ exists and is bounded on $\vL^{2,-s}$.
Since 
\be
 1+ R_0(z) W - \f{1}{\sqrt z} G_{-1}W =1 + G_0W + O(|z|^{\ep})
\ee
  in $\vL^{2,-s}$ for $z \in U_\delta$,  $1+ R_0(z) W - \f{1}{\sqrt z} G_{-1}W$ is invertible for $z\in U_\delta$ if $\delta>0$ is small enough.  
  Denote
\be
D(z) =\left(1+ R_0(z) W - \f{1}{\sqrt z} G_{-1}W\right)^{-1}.
\ee
If $\rho > 2k +2$, one has
\be \label{2.5.20}
D(z) = D_0  + \sum_{j=1}^k z^{\f j 2} D_j + O(|z|^{k+\ep})
\ee
in $\vB(1,-s; 1,-s)$ for  $ k +\f 3 2 <s < \f {\rho+1} 2$,  where
\bea
D_0 &= & (1+G_0W)^{-1} \\
D_1 &= & - D_0 G_1 WD_0 \\
D_2 &=&  (D_0G_1W)^2D_0  - D_0G_2 WD_0
\eea
It follows that 
\be
(1+ R_0(z)W )^{-1}= D(z)( 1 +M(z))^{-1}
\ee
where $M(z) = \f{1}{\sqrt z} G_{-1}WD(z)$. $M(z)$ is  an operator of  rank one.  In order to study the invertibility of $1 + M(z)$, consider the equation
\be \label{2.5.21}
(1+M(z)) u = f,
\ee
where $f \in \vL^{2,-s}$ is given and $u =u(z)$ is to be determined. Take $\phi^*(x, v) =\chi(x) \psi_0(v) $ with $\chi\in \vS(\bR)$ such that 
\[
\int_\bR \chi(x) dx =1.
\] 
Let  $\Pi_0 =  \w{\cdot, \phi^*}\gm_0$. Then $\Pi_0^2 =\Pi_0$. 
Decompose $f$ and $u$ as $f = f_0 + f_1$ and $ u= u_0 + u_1$ where $f_0 = \Pi_0f$, $f_1 = (1-\Pi_0)f$,  and similarly for $u$.
Equation (\ref{2.5.21}) is equivalent with
\be \label{2.5.22}
u_1 =  f_1 \mbox{ and  }  
\ee
\be \label{2.5.22b}
 C(z) (1 + \w{M(z)\gm_0, \phi^*}) = \w{f, \phi^*} -\w{M(z)f_1, \phi^*}
\ee
where $C(z) =\w{u,\phi^*}$ is some constant to be calculated. If  $1 + \w{M(z)\gm_0, \phi^*}\neq 0$ for $z\in U_\delta$ , as we shall prove below, then $C(z)$ is uniquely determined by (\ref{2.5.22b}). Consequently, the equation $(1+ M(z))u =f$ has a unique solution given by
\be
u = C(z) \gm_0 + f_1.
\ee
This will prove the invertibility of $1+ M(z)$ for $z\in U_\delta$. 
\\

Let  us now study  
\be
m(z)= 1 + \w{M(z)\gm_0, \phi^*}
\ee
for $z\in U_\delta$.  Applying (\ref{2.5.20}) with $k=1$ (we need here the condition $\rho>4$), one obtains
\be
\w{M(z)\gm_0, \phi^*} = \f i{2\sqrt z} \w{W D(z)\gm_0, \gm_0} = \f i{2\sqrt z} \left(\sigma_0 + \sqrt z \sigma_1 +  O(|z|^{\f 1 2 +\ep})\right)
\ee
where $\sigma_j =  \w{W D_j\gm_0, \gm_0}$.   By Theorem \ref{th3.4.7}, 
\be
 (1+G_0W)^{-1}\gm_0 =\gm.
\ee 
 Consequently 
\be 
\sigma_0 = \w{W\gm, \gm_0}=0 
\ee
and 
\beas
\sigma_1 & = &   \w{(1+G_0W)^{-1} G_1 W(1+G_0W)^{-1}\gm_0, W\gm_0} \nonumber \\
&=&  \w{G_1 W \gm,  D_0^*W\gm_0}
\eeas
Let $J$ be the symmetry in velocity variable defined by $J : g(x, v) \to (Jg)(x, v)=  g(x,-v)$. Then $J^2 =1$ and 
\be
JPJ =P^*, \quad JWJ = - W \quad  \mbox{ and } JP_0J =P_0^*.
\ee
It follows that $(R_0(z)W)^* = J W R_0(\overline z) J$, hence
\be
(1+ G_0 W)^* = J(1+ WG_0)J.
\ee
We derive that
\beas
D_0^*W\gm_0 &=  & J (1+ WG_0)^{-1}JW\gm_0 \\
 &= &  -J(1+WG_0)^{-1}W\gm_0  \\
 & = &  -J W(1+G_0W)^{-1}\gm_0 = -J W \gm = W\gm.
\eeas
This shows
\be
\sigma_1 = \w{G_1W\gm, W\gm}.
\ee
Since $G_1 =\f 1{\s z} (R_0 (z) - \f 1{\s z}G_{-1} -G_0) + O(|z|^\ep)$ in $\vB(-1, s; 1,- s)$, $s>\f 5 2$,  noticing that
$G_{-1}W\gm=0$, $(1+G_0W)\gm =\gm_0$,  one obtains for $z =\la <0$
\bea
\w{G_1W\gm, W\gm} &= & -\f i{\s{|\la|}} \w{R_0(\la)W\gm, W\gm} + O(|\la|^\ep) \\
 & = &  i\s{|\la|} \w{R_0(\la)\gm, W\gm} + O(|\la|^\ep).
\eea

\begin{prop}\label{prop3.4.8} Assume $\rho>4$. One has
\be
\w{G_1W\gm, W\gm} = i\lim_{\la \to 0_-}  \s{|\la|} \w{R_0(\la)\gm, W\gm} =0.
\ee
\end{prop}
\pf Let $\la <0$ and $\Pi'$  be defined as in the proof of Theorem \ref{th3.4.7}. Then
$\w{R_0(\la)\gm, W\gm}= \w{\Pi'R_0(\la)\gm, W\gm}$, since $\w{\gm_0, W\gm}=0$.
 One has
\be
 R_0(\la)\gm = ( b_0^w(v,D_x,D_v) (-\Delta_x- \lambda)^{-1} + r_0(\lambda))\gm
\ee
in $\vL^{2,-r}$ for any $r>\f 1 2$ and it follows from (\ref{eq3.4.51}) that
\be
\s{|\la|} r_0(\la)\gm = O(\s{|\la|})
\ee
in $\vH^{1,-r}$.  Let us evaluate $\s{|\lambda|}\Pi'b_0^w(v,D_x,D_v) (-\Delta_x- \lambda)^{-1} \gm$. 
\beas
\lefteqn{\s{|\lambda|}  \Pi'b_0^w(v,D_x,D_v) (-\Delta_x- \lambda)^{-1} \gm (x,v)}  \\
& =& 
\f{i} 2  \int_{\bR^4}( e^{-\sqrt{|\lambda|}|y-(x-x')|} -  e^{-\sqrt{|\lambda|}|y-(y'-x')|}) \phi(y')\Phi(v,v',y)\gm(x',v')\; dy dy'dx' dv' 
 \\
 & = &  \f{i} 2  \int_{\bR^4}( e^{-\sqrt{|\lambda|}|y-(x-x')|} -  e^{-\sqrt{|\lambda|}|y-(y'-x')|}) \phi(y')\Phi(v,v',y)\gm_0(v')\; dy dy'dx' dv' \\
  & +  & \f{i} 2  \int_{\bR^4}( e^{-\sqrt{|\lambda|}|y-(x-x')|} -  e^{-\sqrt{|\lambda|}|y-(y'-x')|}) \phi(y')\Phi(v,v',y)(\gm(x',v)-\gm_0(v'))\; dy dy'dx'  \\ 
 &= & \f{i} 2  \int_{\bR^4}( e^{-\sqrt{|\lambda|}|y-(x-x')|} -  e^{-\sqrt{|\lambda|}|y-(y'-x')|}) \phi(y')\Phi(v,v',y)(\gm(x',v)-\gm_0(v'))\; dy dy'dx' \\
 &= &  O(\s{|\la|}|x|\psi_0(v))
\eeas
for $(x, v)\in \bR^2$. The first term on the right-hand side of the second equality above  vanishes by first integrating  with respect to $x'$ variable.
 In the last equality above, we used the upper bound
\[
 |e^{-\sqrt{|\lambda|}|y-(x-x')|} -  e^{-\sqrt{|\lambda|}|y-(y'-x')|}| \le \s{|\la|}|x-y'| \left( e^{-\sqrt{|\lambda|}|y-(x-x')|} +  e^{-\sqrt{|\lambda|}|y-(y'-x')|} \right)
\]
and the fact $\gm-\gm_0 =  O(\w{x}^{-\rho})\psi_0(v)$  to evaluate the integral.
It follows that
\[
\s{|\la|}\w{\Pi'R_0(\la)\gm, W\gm} = O(\s{|\la|}), \quad \la \to 0_-
\]
which finishes the proof of Proposition \ref{prop3.4.8}.
\ef

Summing up, we proved that if $\rho>4$,  then  $m(z)= 1 +  \f i{2\sqrt z} \w{W D(z)\gm_0, \gm_0}$ verifies
\be \label{mz}
m(z)= 1 +   O( |z|^\ep), \quad \ep>0,
\ee
for $z\in U_\delta$.  Therefore $1+M(z)$ is invertible for $z\in U_\delta$  with $\delta>0$  small enough and the solution $u$ to the equation  $(1+M(z))u =f$ is given by
\bea
u& =&  f_1 +   \f 1{m(z)} ( \w{f, \phi^*} -\w{M(z)f_1, \phi^*})\vp_0 \nonumber \\
 &= &  f - \w{f, \phi^*}\gm_0  +   \f 1{m(z)} ( \w{f, \phi^*} -\w{M(z)(f- \w{f, \phi^*}\gm_0), \phi^*})\gm_0 \nonumber \\
 & = & f -  \f 1{m(z)} \w{M(z)f, \phi^*} \gm_0.
\eea
Taking notice that $\w{\gm_0,\phi^*} =1$, we proved the following \\

\begin{prop}\label{prop3.4.9} Let $\rho>4$. Then $1+M(z)$ is invertible in $\vB(1,-s; -1, s)$, $s>\f 3 2$, for $z\in U_\delta$. Its inverse is given by
\be
(1+ M(z))^{-1} = 1 - \f 1{m(z)\s z}  G_{-1} W D(z). 
\ee
 In addition, if $\rho > 2k +2$ for some $k\ge 1$, one has
\be
(1+ M(z))^{-1} = 1 -  \f 1{m(z)\s z}  G_{-1} W \left(D_0  + \sum_{j=1}^k z^{\f j 2} D_j + O(|z|^{k+\ep})\right)
\ee
in $\vB(1,-s; 1,-s)$ for  $ k +\f 3 2 <s < \f {\rho + 1} 2$, where $D_j$ is given by \eqref{2.5.20}.
\end{prop}

\begin{theorem}\label{th3.4.10} Let $\rho>4$. Then there exists some constant $ \delta>0 $ such that
 if $s > \f 5 2$ 
 \be \label{3.4.89}
R(z) = \f{i}{2\sqrt z} \w{\cdot, \gm}\gm + O(|z|^{-\f 1 2 +\ep}),  \quad z\in U_\delta,
\ee
in $\vB(-1,s; 1,-s)$  for some $\ep>0$. In particular, $P$ has no eigenvalues in $U_\delta$. In addition,  the boundary values $R(\la \pm i0)$  of $R(z)$  exist in $\vB(-1, s; 1,-s)$, $s >\f 3 2$, for $\la \in ]0, \delta[$ and is H\"older continuous in $\la \in ]0, \delta[$.
\end{theorem}
\pf We see from the above calculation that $(1+ M(z))^{-1}$ admits an asymptotic expansion as $z\in U_\delta$ and $z\to 0$. 
The existence of  the asymptotics of the resolvent $R(z)$ follows from the equation
\be \label{e5.26}
R(z) = D(z)( 1 +M(z))^{-1}R_0(z)= D(z)\left(1 - \f 1{m(z)\s z}  G_{-1} W D(z) \right)R_0(z).
\ee
Let us calculate its leading term.
\beas
\lefteqn{\left(1 - \f 1{m(z)\s z}  G_{-1} W D(z)\right )R_0(z)}\\
&\equiv &  -\f 1 {m(z)z}  G_{-1}WD_0G_{-1}  + \f 1 {\s z} \left(G_{-1} - \f 1{m(z)} ( G_{-1}WD_0 G_0 +  G_{-1}WD_1 G_{-1})\right).
\eeas
Here and in the following,  ``$\equiv$" means equality module some term which is of order $O(|z|^{-\f 1 2 +\ep})$ in $\vB(-1,s; 1,-s)$, $s>\f 5 2$.
Recall that $G_{-1} =\f i 2 \w{\cdot, \gm_0} \gm_0$, $D_0 = (1+G_0W)^{-1}$ and $(1+G_0W)^{-1}\gm_0 =\gm$. It follows that
\be
G_{-1}WD_0G_{-1} =\f i 2 \w{W\gm, \gm_0} \w{\cdot, \gm_0}\gm_0  =0.
\ee
Consequently 
\beas
\lefteqn{D(z)\left(1 - \f 1{m(z)\s z}  G_{-1} W D(z) \right)R_0(z)}\\
& \equiv &  \f 1{\s z}  D_0 \left(G_{-1} - \f 1{m(z)} ( G_{-1}WD_0 G_0 +  G_{-1}WD_1 G_{-1})\right)
\eeas
Noticing that $m(z) = 1 + O(|z|^{\ep})$, one obtains
\bea \label{3.4.41}
R(z)
  &\equiv &
 \f 1{\s z}  D_0 G_{-1} ( 1 - W( D_0 G_0 +  D_1 G_{-1}))\\
 & = & \f i {2\s z} \w{(( 1 - W( D_0 G_0 +  D_1 G_{-1}))\cdot, \gm_0}\gm \nonumber
\eea
Recall that $D_0^* W\gm_0 = W\gm$ and $\w{ G_1W\gm, W\gm}=0$ (see Proposition \ref{prop3.4.8}). 
One can simplify the leading term as follows:
\beas
\lefteqn{\w{( 1 - W( D_0 G_0 ))\cdot, \gm_0}}
\\
 & = & \w{ \cdot, \gm_0} +   \w{ \cdot, G_0^* D_0^* W\gm_0} 
 =  \w{ \cdot, \gm_0} +   \w{ \cdot, G_0^* W\gm} \\
& = & \w{ \cdot, \gm_0} +   \w{ \cdot, JG_0J W\gm} 
 =  \w{ \cdot, \gm_0} -   \w{ \cdot, G_0 W\gm} =  \w{\cdot, \gm}
\eeas
and 
\beas
\lefteqn{\w{ WD_1 G_{-1}\cdot, \gm_0}}\\
 &= & - \w{ W D_0G_1WD_0 G_{-1}\cdot, \gm_0} 
  =  -\f i 2 \w{\cdot, \gm_0} \w{ W D_0G_1W\gm, \gm_0} \\ 
 & = &  \f i 2 \w{\cdot, \gm_0} \w{ G_1W\gm, D_0^* W\gm_0} 
   =  \f i 2 \w{\cdot, \gm_0} \w{ G_1W\gm, W\gm} =0.
\eeas
This finishes the proof of \eqref{3.4.89}. \eqref{3.4.89} implies that $R(z)$ has no poles in $U_\delta$, hence $P$ has no eigenvalues there. The last statement of Theorem \ref{th3.4.10} is a consequence of Corollary \ref{cor2.2} (b) and (\ref{e5.26}), since the boundary values $D(\la \pm i0)$ exist in $\vB(1, -s; 1,-s)$, $s >\f 3 2$, for $\la \in ]0, \delta[$ and are continuous in $\la$.
\ef

\sect{Large time asymptotics of solutions} 

The following high energy resolvent estimate is proved in \cite{w3}. 

\begin{theorem}\label{th5.1} Let $n \ge 1$ and assume (\ref{ass1}) with $\rho \ge -1$. Then there exists $C>0$ such that   $\sigma (P) \cap \{z; |\Im z| > C,  \Re z \le \f 1 C |\Im z|^{\f 1 2}  \} = \emptyset$ and
\be \label{eq5.1}
\|R(z)\| \le \f{C}{|z|^{\f 1 2}}, \quad
\ee
and
\be \label{eq5.2}
\|(1-\Delta_v + v^2)^{\f 1 2}R(z)\| \le \f{C}{|z|^{\f 1 4}}, \quad
\ee
for $|\Im z| > C $ and $ \Re z \le \f 1 C |\Im z|^{\f 1 2}$.
\end{theorem}

Let $S(t) =e^{- t P}$, $t\ge 0$, be the one-parameter strongly continuous semigroup generated by $-P$. Then one can firstly represent $S(t)$ as
\be \label{eq5.3}
S(t) f =  \frac{1}{2\pi i} \int_{\gamma } e^{-t z}R(z)f dz
\ee
for $f\in L^2(\bR^2)$ and $t>0$, where  the contour $\gamma$ is chosen such that
\[\gamma = \gamma_- \cup \gamma_0 \cup \gamma_+\]
with $\gamma_{\pm} =\{ z; z = \pm i C + \lambda \pm i C \lambda^2, \lambda \ge 0\}$ and $\gamma_0$ is a curve in the left-half complexe plane joining $-i C$ and $iC$ for some $C>0$ sufficiently large, $\gamma$ being oriented from $-i\infty$ to $+i\infty$.

Remark that under the condition (\ref{ass1}) with $\rho>0$, $P$ has no eigenvalue on the imaginary axis (\cite{hln}). 
Making use of  analytic deformation  and Theorem \ref{th3.4.10}, one obtains from (\ref{eq5.1}) that
\be \label{6.4}
\w{S(t)f,g} =  \frac{1}{2\pi i} \int_{\Gamma} e^{-t z}\w{R(z)f,g} dz, \quad t>0,
\ee
for any $f,g \in \vL^{2, s}$ with $s> \f 5 2$. Here 
\[
\Gamma = \Gamma_- \cup \Gamma_0 \cup \Gamma_+
\]
with 
\[
\Gamma_{\pm} =\{ z; z = \delta + \lambda \pm i  \delta^{-1} \lambda^2, \lambda \ge 0\}
\]
 for $\delta>0$ small enough and 
 \[
 \Gamma_0 =\{ z = \lambda + i0; \lambda \in [0, \delta]\}\cup \{ z = \lambda - i0; \lambda \in [0, \delta]\}.
 \] 
 $\Gamma$ is oriented from $-i\infty$ to $+i\infty$.\\

{\noindent \bf Proof of Theorem \ref{th1.1}.}
By \eqref{6.4}, one has  for $f,g \in \vL^{2,s}(\bR^2)$ with $s > \frac{5}{2}$
\beas
   \w{S(t)f, g} &=& \frac{1}{2 \pi i} \left(\int_{\Gamma_0} + \int_{\Gamma_-} + \int_{\Gamma_+} \right)  e^{-tz} \w{R(z) f, g} \, d z \\
   & \equiv & I_1 + I_2 + I_3.
\eeas
For $\delta>0$ appropriately small and fixed, it follows  from Theorem \ref{th5.1} that there exist some constants $C, c>0$ such that
\be
|I_j| \le C e^{-c t} \|f\|\;  \|g\|, \quad t >0,
\ee
 for $j=2,3$.  Set
 \[
 F_{-1} = \f i 2 \w{\cdot, \gm}\gm.
 \]
 Applying Theorem \ref{th3.4.10}, one has
  \begin{equation}
    \begin{aligned}
    I_1 &= \frac{1}{2\pi i} \int_0^{\delta} e ^{-t \lambda} \w{(R(\la+ i0)-R(\la -i0))f, g} \, d \lambda\\
    & = \frac{1}{\pi i} \int_0^{\delta} e ^{-t \lambda} \la^{-\f 1 2}\w{(F_{-1} + O(\la^\ep))f, g} \, d \lambda
     \\ 
     & = \frac{1}{\pi i }\int_0^{+\infty} \frac{1}{\sqrt{\lambda}}e ^{-t \lambda} \w{F_{-1}f,  g} \, d \lambda +  O(t^{-\frac{1}{2} - \ep})\|f\|_{0,s}\|g\|_{0,s}\\
    &= \frac{1}{i\sqrt{\pi t}  } \w{F_{-1}f, g} + O(t^{-\frac{1}{2} - \ep})
    \end{aligned}
    \end{equation}
    as $t \rightarrow +\infty$ for some $\ep > 0$. Using the formula for $F_{-1}$, we arrive at
    \begin{equation}
      S(t) =  \frac{1}{(4\pi t)^{\frac{1}{2}}} \w{\cdot, \gm} \gm + O(t^{-\frac{1}{2} - \ep}), t \to +\infty
    \end{equation}
as operators in $\vB(0,s; 0,-s)$ with $s>\f 5 2$. Theorem \ref{th1.1} is proved.
\ef

\end{document}